\documentclass[reqno,12pt]{amsart}

\usepackage{amsmath}
\usepackage{amssymb}

\renewcommand{\bar}{\overline}

\newcommand{\eps}{\varepsilon}

\newcommand{\CC}{\mathbb{C}}

\newcommand{\FF}{\mathbb{F}}

\newcommand{\HH}{\mathbb{H}}

\newcommand{\PP}{\mathbb{P}}
\newcommand{\QQ}{\mathbb{Q}}
\newcommand{\RR}{\mathbb{R}}

\newcommand{\ZZ}{\mathbb{Z}}

\newcommand{\Qp}{\QQ_p}
\newcommand{\Zp}{\ZZ_p}
\newcommand{\Cp}{\CC_p}

\newcommand{\CK}{\CC_K}
\newcommand{\Qpbar}{\bar{\QQ}_p}
\newcommand{\Fp}{\FF_p}
\newcommand{\Fpbar}{\bar{\FF}_p}

\newcommand{\ints}{{\mathcal O}}

\newcommand{\maxid}{{\mathcal M}}

\newcommand{\khat}{\hat{k}}

\newcommand{\calJ}{{\mathcal J}}

\newcommand{\Fat}{{\mathcal F}}
\newcommand{\Jul}{{\mathcal J}}

\newcommand{\PCK}{\PP^1(\CK)}

\newcommand{\PC}{\PP^1(\CC)}

\newcommand{\Pk}{\PP^1(k)}
\newcommand{\Pkhat}{\PP^1(\hat{k})}

\newcommand{\PGL}{\hbox{\rm PGL}}

\newcommand{\diam}{\hbox{\rm diam}}

\DeclareMathOperator{\charact}{char}

\newcommand{\Dbar}{\bar{D}}

\newtheorem{thm}{Theorem}[section]
\newtheorem{defin}[thm]{Definition}

\newtheorem{lemma}[thm]{Lemma}

\newtheorem*{thmA}{Theorem A}
\newtheorem*{thmB}{Theorem B}

\newcounter{exno}
\newenvironment{example}{\refstepcounter{exno}
	\smallskip {\bf Example \theexno .}}{}

\title[Wandering Domains and Nontrivial Reduction]
{Wandering Domains and Nontrivial Reduction
	in Non-Archimedean Dynamics}
\author{Robert L. Benedetto}
\date{May 12, 2004; revised November 24, 2004}
\address{Department of Mathematics and Computer Science \\
        Amherst College \\
        Amherst, MA 01002 \\
        USA}
\email{rlb@cs.amherst.edu}
\urladdr{http://www.cs.amherst.edu/\textasciitilde rlb}

\begin{document}

\newcounter{bean}
\newcounter{sheep}

\begin{abstract}
Let $K$ be a non-archimedean field with residue field $k$,
and suppose that $k$ is not an algebraic extension of a
finite field.  We prove two results concerning wandering
domains of rational functions $\phi\in K(z)$ and
Rivera-Letelier's notion of nontrivial reduction.
First, if $\phi$ has nontrivial reduction, then assuming some
simple hypotheses, we show that the Fatou set of
$\phi$ has wandering components by any of the usual
definitions of ``components of the Fatou set''.  Second,
we show that if $k$ has characteristic zero and $K$
is discretely valued, then the existence of a wandering
domain implies that some iterate has nontrivial reduction
in some coordinate.
\end{abstract}

\maketitle

The theory of complex dynamics in dimension one,
founded by Fatou and Julia in the early twentieth century,
concerns the action of a rational function $\phi\in\CC(z)$
on the Riemann sphere $\PC=\CC\cup\{\infty\}$.
Any such $\phi$ induces a partition of the
sphere into the closed Julia set $\Jul_{\phi}$,
where small errors become arbitrarily large under
iteration, and the open Fatou set $\Fat_{\phi}=\PC\setminus\Jul_{\phi}$.
There is also a natural action of $\phi$ on the connected
components of $\Fat_{\phi}$, taking a component $U$ to $\phi(U)$,
which is also a connected component of the Fatou set.
In 1985, using quasiconformal methods, Sullivan \cite{Sul}
proved that $\phi\in\CC(z)$ has no wandering domains;
that is, for each component $U$ of $\Fat_{\phi}$, there are integers
$M\geq 0$ and $N\geq 1$ such that $\phi^M(U)=\phi^{M+N}(U)$.
We refer the reader to \cite{Bea,CG,Mil} for
background on complex dynamics.

In the past two decades, there have been a number of
investigations of dynamics over complete metric fields other than
$\RR$ or $\CC$.  All such fields are {\em non-archimedean};
that is, the metric on the field $K$
satisfies the {\em ultrametric} triangle inequality
$$d(x,z) \leq \max \{ d(x,y) , d(y,z) \}
\quad \text{for all } x,y,z\in K.$$
Herman and Yoccoz \cite{HY} first considered dynamics over
such fields in a study of
linearization at fixed points, in part to discover
which properties of complex dynamical systems are specific
to archimedean fields and which are more general.
The question of comparing archimedean and non-archimedean
dynamics has continued to drive the field, as have questions
arising in number theory in the study of rational dynamics;
\cite{Ben3,Ben4,Ben5,Bez,Hs2,MorSil2,Riv1,Riv2,SW}.

In particular, it is natural to ask how the dynamical
properties of Fatou components
extend to the non-archimedean setting.  In \cite{Ben2}, the author
proved a no wandering domains theorem over $p$-adic
fields, assuming some weak hypotheses.  That theorem
relied heavily on the fact that the residue field $k$
(see below) of
a $p$-adic field $K$ is an algebraic extension of the
finite field $\FF_p$.  In fact, for non-archimedean fields $K$
without such a residue field, it is easy to construct
rational functions with wandering domains; see 
Example~\ref{ex:fixhpwd} and \cite[Example~2]{Ben5}.
The aim of this paper is to classify all such wandering
domains.

We fix the following notation.
\begin{tabbing}
$K$ \= \hspace{1.0in} \=
  a complete non-archimedean field with absolute value $|\cdot|$ \\
$\hat{K}$ \> \> an algebraic closure of $K$ \\
$\CK$ \> \> the completion of $\hat{K}$ \\
$\ints_K$ \> \> the ring of integers $\{x\in K : |x|\leq 1\}$ of $K$ \\
$k$ \> \> the residue field of $K$ \\
$\ints_{\CK}$ \> \> the ring of integers $\{x\in \CK : |x|\leq 1\}$ of $\CK$ \\
$\hat{k}$ \> \> the residue field of $\CK$ \\
$\PCK$ \> \> the projective line $\CK\cup\{\infty\}$
\end{tabbing}
Recall that the absolute value $|\cdot|$ extends in
unique fashion to $\hat{K}$ and to $\CK$.
Recall also that the residue field $k$ is defined to
be $\ints_K/\maxid_K$, where
$\maxid_K$ is the maximal ideal $\{x\in K : |x|< 1\}$ of $\ints_K$.
The residue field $\hat{k}$ is defined similarly.
There is a natural inclusion of the residue field $k$ into
$\hat{k}$, making $\hat{k}$ an algebraic closure of $k$.
We refer the reader to \cite{Gou,Kob,Rob,Ser2}
for surveys of non-archimedean fields.

The best known complete non-archimedean field is
$K=\Qp$, the field of $p$-adic rational numbers
(for any fixed prime number $p$).
Its algebraic closure is $\hat{K}=\Qpbar$, and the completion
$\CK$ is frequently denoted $\Cp$.
The ring of integers is $\ints_K=\Zp$, with residue
field $k=\FF_p$ (the field of $p$ elements), and $\hat{k}$ is
the algebraic closure $\Fpbar$.
Note that $\charact K=0$, but $\charact k=p$.  Thus, we say
the characteristic of $\Qp$ is $0$, but the {\em residue characteristic}
of $\Qp$ is $p$.

As another example,
if $L$ is any abstract field, then $K=L((T))$, the field of formal
Laurent series with coefficients in $L$, is a complete non-archimedean
field with $\ints_K=L[[T]]$ (the ring of formal Taylor series)
and $k=L$.  In this case, $\charact K=\charact k = \charact L$.
The absolute value $|\cdot|$ on $K$ may be
defined by $|f| = 2^{-n}$, where $n\in\ZZ$ is the least
integer for which the $T^n$ term of the formal Laurent
series $f$ has a nonzero coefficient.
Dynamics over such a function field $K$ has
applications to the study of one-parameter families
of functions defined over the original field $L$; see,
for example, \cite{Kiwi}.

In the study of non-archimedean dynamics of
one-variable rational functions,
we consider a rational function
$\phi\in K(z)$, which acts on $\PCK$
in the same way that a complex rational function
acts on the Riemann sphere.
In \cite{Ben1,Ben2,Ben5},
the author defined non-archimedean Fatou and Julia sets.
The wandering
domains we will study are components of the Fatou set,
which is an open subset of $\PCK$.

Of course, to study wandering domains, we must first have
an appropriate notion of ``connected components'' of subsets
of $\PCK$, which is a totally disconnected topological space.
Several definitions have been proposed in the literature,
and each is useful in slightly different settings, just as
connected components and path-connected components are related
but distinct notions.  We will consider four definitions
in this paper, all of which are closely related and which
frequently coincide with one another.

In \cite{Ben1,Ben2,Ben5}, the author
proposed two analogues of ``connected components''
of the Fatou set: D-components and analytic components,
both of which we will define in Section~\ref{sect:back}.
Rivera-Letelier proposed an alternate
definition in \cite{Riv1,Riv2}.  His definition
was stated only over the $p$-adic field $\Cp$;
in Section~\ref{sect:back},
we define an equivalent version of his components,
which we call dynamical components, for all
non-archimedean fields.
We also will propose a fourth analogue,
called dynamical D-components, which will
actually be useful for proving things
about the other three types of components.

Our first main result
generalizes the aforementioned
wandering domain from \cite[Example~2]{Ben5},
for the function $\phi(z) = (z^3 + (1+T)z^2)/(z+1)\in K(z)$,
where $K=\QQ((T))$.
In that example,
the wandering domain $U$ and
all of its foward images $\phi^n(U)$
are open disks of the form
$D(a,1)$, with $|a|\leq 1$, where $D(a,r)$ denotes the
open disk of radius $r$ about $a$.  In fact, the map $\phi$
has the property that for all but
finitely many of the disks $D(a,1)$ with $|a|\leq 1$,
the image $\phi(D(a,1))$ is just $D(\phi(a),1)$;
in Rivera-Letelier's language \cite{Riv1,Riv2},
$\phi$ has {\em nontrivial reduction}.
Equivalently, $\phi$ has a fixed point in Rivera-Letelier's
``hyperbolic space'' $\HH$, which is essentially the same
space as the Berkovich projective line $\PP^1_{\text{Berk}}(\CK)$;
see \cite{Ber,RuBa}.
Both $\HH$ and $\PP^1_{\text{Berk}}(\CK)$ 
have been used with increasing frequency in the study
of the mapping properties and dynamics of non-archimedean
rational functions.  Readers familiar with either space
will recognize the wandering domains we will construct as
connected components of the full space $\PP^1_{\text{Berk}}(\CK)$
with a single type~II point removed.  However, for
simplicity, we will present
our arguments without reference to Berkovich spaces.

We will define nontrivial reduction precisely in
Section~\ref{sect:nontriv}.
We will then present Theorem~\ref{thm:fixhpwd}
and Example~\ref{ex:fixhpwd}, which imply
the following result.
\begin{thmA}
Let $K$ be a non-archimedean field
with residue field $k$, where $k$
is not an algebraic extension of
a finite field.
\begin{list}{\rm (\alph{bean}).}{\usecounter{bean}}
\item
Let $\phi\in K(z)$ be a rational function
of nontrivial reduction $\bar{\phi}$ with
$\deg\bar{\phi}\geq 2$.
Then $\phi$ has a wandering dynamical component $U$
such that for every integer $n\geq 0$, $\phi^n(U)$
is an open disk of the form $D(b_n,1)$, with $|b_n|\leq 1$.
$U$ is also a wandering dynamical D-component;
moreover, if the Julia set of $\phi$ intersects
infinitely many residue classes of $\PCK$, then
$U$ is also
a wandering D-component, and a wandering analytic component.
\item
There exist functions $\phi\in K(z)$ satisfying all
of the hypotheses of part~{\rm (a)}.
\end{list}
\end{thmA}
Theorem~\ref{thm:fixhpwd} is an even stronger result:
that under the hypotheses of Theorem~A, there
are actually infinitely many different grand orbits of
wandering domains of the form $D(b,1)$.
In addition, Theorem~\ref{thm:julia} will give
sufficient conditions for the Julia set of $\phi$
to intersect infinitely many residue classes.

Theorem~A and Theorem~\ref{thm:fixhpwd} apply to
maps with reduction $\bar{\phi}$ of degree at least two.
If $\phi$ has a nontrivial reduction $\bar{\phi}$ of degree one,
the situation is a bit more complicated.
Examples~\ref{ex:degonewd}--\ref{ex:degonestat},
will show that in some such cases there is a wandering
domain of the form $D(b,1)$, and in other cases
there is not.

Given Theorem~A, we can construct still more wandering
domains over appropriate fields $K$, as follows.
If $\phi\in K(z)$ is a rational function
and $g\in\PGL(2,\CK)$ is a linear fractional transformation,
suppose that some
conjugated iterate $\psi(z)=g\circ \phi^n\circ g^{-1}(z)$
has nontrivial reduction of degree at least two.
Even if the original map $\phi$
has trivial reduction, Theorem~\ref{thm:fixhpwd}
shows that $\phi$ has a wandering domain, because $\psi$ does.

The existence of rational functions with
wandering domains should not come as
a surprise for fields $K$ satisfying the
hypotheses of Theorem~A.  For example, the infinite
residue field prevents $K$ from being locally
compact, allowing plenty of room for the various iterates
of the wandering domain to coexist.
Thus, the impact of Theorem~A
is not so much the fact that wandering domains exist, but that
they may be produced by such simple reduction
conditions.

Perhaps more interesting than
the existence of such wandering domains is our next
theorem, which shows that for certain fields $K$,
the {\em only} wandering domains possible
for rational functions are those
described above.  That is, any wandering domain
for such a field {\em must} come from a nontrivial
reduction.

\begin{thmB}
Let $K$ be a non-archimedean field with residue field $k$.
Suppose that $K$ is discretely valued and that $\charact k = 0$.
Let $\phi(z)\in K(z)$ be a rational function, and suppose
that some $U\subset\PCK$ is a wandering domain
(analytic, dynamical, D-component, or dynamical D-component)
of $\phi$.  Then there are
integers $M\geq 0$ and $N\geq 1$ and a change of coordinate
$g\in\PGL(2,\CK)$ with the following property:

Let $\psi(z)=g\circ\phi^N\circ g^{-1}(z)\in\CK(z)$.
Then $D(0,1)$ is the component of the Fatou set of $\psi$
containing $g(\phi^M(U))$, and $\psi$ has nontrivial reduction.
\end{thmB}

The clause ``$D(0,1)$ is the component of the Fatou set of $\psi$
containing $g(\phi^M(U))$'' is equivalent to ``$g(\phi^M(U))=D(0,1)$''
if we are dealing with analytic or dynamical components.  On the
other hand, if $U$ is a D-component or dynamical D-component, then
it is possible that $\phi^M(U)$ is a proper subset of a
(dynamical) D-component.
We will prove a slightly stronger version of Theorem~B
in Theorem~\ref{thm:nowd} of Section~\ref{sect:czero},
showing the the function $g$ can be defined over a finite
extension of $K$.

In \cite{Ben6,Ben8}, it was shown that
rational functions, including polynomials,
may have wandering domains if
$K=\CK$ and if $\charact k>0$.
(The hypotheses of the no wandering domains
theorems of \cite{Ben1,Ben2,Ben5} require
the field of definition
$K$ to be locally compact, whereas $\CK$
is not locally compact.)
By contrast, Theorem~B shows that
rational functions have no wandering domains (besides those
arising from a nontrivial reduction, as in
Example~\ref{ex:fixhpwd}) if the field of
definition $K$ is discretely valued and has
residue characteristic {\em zero}.
However, an algebraically closed non-archimedean
field (such as $\CK$) cannot be discretely valued.
Thus, one is naturally
led to ask the following open questions:
\begin{enumerate}
\item If $K$ is locally compact (implying both
that $K$ is discretely valued and has residue
characteristic $p>0$), do there exist polynomial
or rational functions
$\phi\in K(z)$ with wandering domains?
\item If $K$ is complete and {\em algebraically closed},
with residue characteristic zero,
do there exist functions $\phi\in K(z)$
with wandering domains other than those
arising from a nontrivial reduction?
\end{enumerate}
The results of \cite{Ben1,Ben2,Ben5} suggest 
that the answer to the first question
above is probably ``no''.  Meanwhile,
we know of no progress on the second question.

The outline of the paper is as follows.
We will begin in Section~\ref{sect:back} with a review
of some dynamical terminology and
the definitions of the various types of Fatou components.
In Section~\ref{sect:nontriv}, we will introduce
Rivera-Letelier's notion of nontrivial reduction
and state three important lemmas.
In Section~\ref{sect:height}, we will recall
a few facts from the theory of diophantine
height functions.  Heights will be used only
in the proof of Lemma~\ref{lem:wandpt};
the reader unfamiliar with the theory may
skip both Section~\ref{sect:height}
and the proof of Lemma~\ref{lem:wandpt} without
loss of continuity.
Finally, in Section~\ref{sect:exist} we will
prove Theorem~A, and
in Section~\ref{sect:czero} we will prove
Theorem~B.
We also include an appendix on the relevant
terminology and fundamental properties of
rational functions and some non-archimedean analysis.

The author would like to thank the referee for
a careful reading of the paper and many helpful
suggestions.

\section{Dynamical Terminology and Fatou Components}
\label{sect:back}

Let $X$ be a set, and let $f:X\rightarrow X$ be a function.
For any $n\geq 1$, we write
$f^1=f$, $f^2=f\circ f$, and in general, $f^{n+1}=f\circ f^n$;
we also define $f^0$ to be the identity function on $X$.
Let $x\in X$.  We say that $x$ is {\em fixed}
if $f(x)=x$; that $x$ is {\em periodic of period $n\geq 1$}
if $f^n(x)=x$; that $x$ is {\em preperiodic}
if $f^m(x)$ is periodic for some $m\geq 0$; or
that $x$ is {\em wandering} if $x$ is not preperiodic.
Note that all fixed points are periodic, and all periodic
points are preperiodic.  We define the {\em forward orbit}
of $x$ to be the set $\{f^n(x) : n\geq 0\}$;
the {\em backwards orbit} of $x$ to be
$\bigcup_{n\geq 0} f^{-n}(x)$; and the {\em grand orbit}
of $x$ to be
$$\{y\in X : \exists\, m,n\geq 0 \text{ such that } f^m(x)=f^n(y)\}.$$
Equivalently, the grand orbit of $x$ is the union of the
backwards orbits of all points in the forward orbit of $x$.
We say a grand orbit $S$ is {\em preperiodic} if it contains
a preperiodic point, or $S$ is {\em wandering} otherwise.
Note that $S$ is preperiodic (respectively, wandering)
if and only if every
point in $S$ is preperiodic (respectively, wandering).

Suppose $X,Y$ are metric spaces.  Recall that the family $F$
of functions from $X$ to $Y$ is {\em equicontinuous} at $x\in X$ if 
for every $\eps>0$ there is a $\delta>0$ such that
$d(f^n(x),f^n(x'))<\eps$ for all $f\in F$ and for all $x'\in X$
satisfying $d(x,x')<\delta$.  (The key point is that the
choice of $\eps$ is independent of $f$.)

Now consider $X=\PCK$ and $f=\phi\in\CK(z)$.  The {\em Fatou set}
of $\phi$ is the set $\Fat=\Fat_{\phi}$ consisting of all
points $x\in\PCK$ for which $\{\phi^n : n\geq 0\}$ is equicontinuous
on some neighborhood of $x$, with respect to the spherical metric
(see, for example, \cite[Section~5]{MorSil2})
on $\PCK$.  The {\em Julia set} $\Jul=\Jul_{\phi}$
of $\phi$ is the complement $\Jul=\PCK\setminus\Fat$.
Clearly the Fatou set is open, and the Julia set is closed.
It is easy to show that $\phi(\Fat)=\phi^{-1}(\Fat)=\Fat$
and $\Fat_{\phi^n}=\Fat_{\phi}$ for all $n\geq 1$, and
similarly for the Julia set.

Intuitively speaking, the Fatou set is the region where small errors
stay small under iteration, while the Julia set is the region of
chaos.  Note that because $\PCK$ is not locally compact, the
Arzel\`{a}-Ascoli theorem fails, which is why non-archimedean
Fatou and Julia sets are defined in terms of equicontinuity
instead of normality.

It is easy to verify (using, for example, \cite[Lemma 2.7]{Ben7}
or other well known lemmas on non-archimedean power series)
that if $U\subset K$ is a disk
and if $R>0$ such that for all $n\geq 0$,
$\phi^n(U)$ is a subset of $\CK$ of diameter at most $R$,
then $U\subset\Fat_{\phi}$.
It follows that if $V\subset K$ is any open set
and if $R>0$ such that for all $n\geq 0$,
$\phi^n(V)\subset K$ and $\diam(\phi^n(V))\leq R$,
then $V\subset\Fat_{\phi}$.
Another criterion, due to Hsia \cite{Hs2} (see
also \cite[Theorem 3.7]{Ben7}) states that
if $U\subset\PCK$
is a disk such that $\bigcup_{n\geq 0} \phi^n(U)$
omits at least two points of $\PCK$, then $U\subset\Fat_{\phi}$.
Clearly Hsia's criterion also extends to arbitrary open sets $V$
in place of $U$.

Using the language of affinoids from Section~\ref{ssec:rigid}
of the Appendix, we now define components of non-archimedean
Fatou sets.
\begin{defin}
\label{def:comps}
Let $\phi\in\CK(z)$ be a rational function with Fatou set $\Fat$,
and let $x\in\Fat$.
\begin{list}{\rm \alph{bean}.}{\usecounter{bean}}
\item The {\em analytic component} of $\Fat$ containing $x$ is the union
  of all connected affinoids $W$ in $\PCK$
  such that $x\in W \subset \Fat$.
\item The {\em D-component} of $\Fat$ containing $x$ is the union
  of all disks $U$ in $\PCK$
  such that $x\in U \subset \Fat$.
\item The {\em dynamical component} of $\Fat$ containing $x$ is the union
  of all rational open connected affinoids $W$ in $\PCK$
  such that $x\in W$ and the set
  $$ \PCK \setminus \left( \bigcup_{n\geq 0} \phi^n(W) \right)$$
  is infinite.
\item The {\em dynamical D-component} of $\Fat$ containing $x$ is the union
  of all rational open disks $U$ in $\PCK$
  such that $x\in U$ and the set
  $$ \PCK \setminus \left( \bigcup_{n\geq 0} \phi^n(U) \right)$$
  is infinite.
\end{list}
\end{defin}
Clearly all of these components are open sets.  Because finite
unions
of overlapping connected affinoids or disks are again connected affionids
or disks (or all of $\PCK$), the relation ``$y$ is in the component
of $\Fat$ containing $x$'' is an equivalence relation between
$x$ and $y$, for each of the four types of components.
Note that by Hsia's criterion, any dynamical component or dynamical
D-component must in fact be contained in the Fatou set, so the
terminology ``component of $\Fat$'' is not misleading.
D-components must be either disks, all of $\PCK$, 
or all but one point of $\PCK$.
Dynamical D-components must be either open disks, all
of $\PCK$, or all but one point of $\PCK$.  Analytic
and dynamical components may be more complicated geometrically.
Frequently, two or more types of components coincide
in a particular case.

Analytic components and D-components were first defined in
\cite{Ben1,Ben2}.  Dynamical components were defined by
Rivera-Letelier in \cite{Riv1,Riv2}; he called them
simply ``components'', and he used a different, but equivalent,
definition.  Dynamical D-components are new to the literature.

For any of the four types of components,
if $x\in\Fat_{\phi}$ and if $V$ is the component containing $x$,
then $\phi(V)$ is contained in the component containing $\phi(x)$,
by of the mapping properties discussed in Sections~\ref{ssec:nonarch}
and~\ref{ssec:rigid} of the appendix.
Thus, $\phi$ induces an action $\Phi_D$ on the
set of D-components by
$$\Phi_D(U) = \text{the D-component containing } \phi(U).$$
Similarly, $\phi$ induces actions
$\Phi_{an}$ on the set of analytic components,
$\Phi_{dyn}$ on the set of dynamical components,
and $\Phi_{dD}$ on the set of dynamical D-components.
Thus, we can discuss fixed components, wandering components,
grand orbits of components, etc., for each of the
four types.

For analytic and dynamical components, it can be shown
\cite{Ben2,Riv1}
that $\Phi_{an}(V)=\phi(V)$ and $\Phi_{dyn}(V)=\phi(V)$.
For D-components and dynamical D-components, the corresponding
equalities usually hold; but occasionally,
the containment may be proper.
Fortunately, by Lemma~\ref{lem:dbound}, for any given $\phi\in\CK(z)$
of degree $d$, there are at most $d-1$ D-components $U$
for which there exists a D-component $V$ with
$\Phi_D(V)=U$ but $\phi(V)\subsetneq U$.
The analogous statement also holds for dynamical D-components.


The following examples should help to clarify how each
of the four types of components behaves.  We omit the proofs
of most of the claims in the following examples.  Details
concerning similar examples may be found in, for example,
\cite{Ben1,Ben2,Ben3,Ben5,Riv1,Riv2}.

\begin{example}
\label{ex:comp1}
Let $n\geq 2$, and let $\phi(z)=z^n$.  Then it is
easy to show that $\Fat_{\phi}=\PCK$.  (The same is
true in the more general situation that $\phi$ has good
reduction; see Section~\ref{sect:nontriv}.)  It follows
immediately that there is only one D-component and only
one analytic component, namely the full set $\PCK$.

On the other hand, all disks of the form $D(\alpha,1)$
for $\alpha\in\CK$ with $|\alpha|\leq 1$, as well as the
disk $\PCK\setminus\Dbar(0,1)$ at $\infty$, 
are dynamical components and dynamical D-components of
the Fatou set.  Indeed, any strictly larger open disk or
affinoid $U$ will have the property that $\bigcup_{n\geq 0}\phi^n(U)$
omits at most the two points $0$ and $\infty$.
(Cf.\ Lemma~\ref{lem:fixclass}.)
\end{example}

\begin{example}
\label{ex:comp2}
Let $p=\charact K\geq 0$, let $c\in K$ with $0<|c|<1$, let $d\geq 2$
be an integer not divisible by $p$,
and let $\phi(z) = z^d - c^{-d}$.  Writing $U_0=\Dbar(0,|c|^{-1})$, it
is easy to see that for $x\in\PCK\setminus U_0$, the iterates
$\phi^n(x)$ approach $\infty$.  Thus, the Julia set $\Jul_{\phi}$
is contained in $U_0$.  In fact, one can check that for any
$n\geq 0$, if we set $U_n = \phi^{-n}(U_0)$, then $U_n$ is
the disjoint union of $d^n$ disks, each of radius $|c|^{n-1}$,
with $d$ such disks in each of the $d^{n-1}$ disks of $U_{n-1}$.
It follows easily that $\Jul_{\phi} =\bigcap_{n\geq 0} U_n$.
The dynamical D-component and the D-component
of $\Fat_{\phi}$ containing $\infty$ are both the
open disk $\PCK\setminus U_0$.  On the other hand, the analytic
and the dynamical component at $\infty$, which also coincide in
this case, are both the more complicated set $\PCK\setminus\Jul_{\phi}$,
which is the whole Fatou set.
\end{example}

\begin{example}
\label{ex:comp3}
Let $p=\charact K$, and assume that $p>0$.
Let $d\geq 2$ be an integer not divisible by $p$,
and let $c\in K$ with $|p|^{1/(pd-1)}<|c|<1$.
Let $\phi(z) = z^{pd}-c^{-pd}$.
Writing $U_0=\Dbar(0,|c|^{-1})$
and $V_0=\PCK\setminus U_0$, it
is again easy to see that for $x\in V_0$, the iterates
$\phi^n(x)$ approach $\infty$.  Defining $V_n=\phi^{-n}(V_0)$ for $n\geq 0$,
the set of points which are attracted to
$\infty$ under iteration is $V=\bigcup_{n\geq 0} V_n$, which is
a complicated union of affinoids (and which is not itself an affinoid).
However, for $z\in U_0$, we have $|\phi'(z)|<1$, from which it
follows easily that $\Jul_{\phi}=\emptyset$.
Thus, the D-component and analytic component of $\Fat_{\phi}$
containing $\infty$ are both $\PCK$ itself, the
dynamical component is $V$, and the dynamical D-component
is $V_0$.
\end{example}

\begin{example}
\label{ex:comp4}
Let $c\in K$ with $0<|c|<1$, 
let $d\geq 3$ be an integer, and let $\phi(z) = cz^d + z^{d-1} + z$.
Clearly $\phi(\Dbar(0,1))=\Dbar(0,1)$.  There is also
a repelling fixed point at $z=-1/c$, and the backwards
orbit of $-1/c$ includes points of aboslute value $|c|^{-1/d^n}$
for arbitrary small $n\geq 0$.  Therefore, any connected
affinoid strictly containing $\Dbar(0,1)$ must intersect
the Julia set.

The dynamical D-component of and the dynamical
component of $\Fat_{\phi}$ containing $0$
are both
the open disk $D(0,1)$.  However, the analytic component
and the D-component of $0$ are both the larger closed disk
$\Dbar(0,1)$.
\end{example}

\begin{example}
\label{ex:comp5}
Let $b,c\in K$ with $0<|c|<|b|=|b-1|=1$, and let
$$\phi(z) = \frac{bz(z+c)(z+c^2)}{(z+bc)(z+c^3)(cz+1)^2}.$$
Then $\phi(0)=0$ is a repelling fixed point, and
$\phi(\infty)=0$.
Let $V$ be the annulus $D(0,|c|^{-1})\setminus \Dbar(0,|c|)$;
then for all $z\in V$, we have $|\phi(z) - bz| < |z|$.
In particular, $\bar{\phi}(z) = \bar{b} z$, and $\phi(V)=V$.
One can also show that there are infinitely
many disks of the form $D(\alpha,|c|)$ and $D(\beta, |c|^{-1})$
with $|\alpha|=|c|$ and $|\beta|=|c|^{-1}$ which contain
preimages of $0$; thus, $\Jul_{\phi}$ intersects infinitely
many such disks.
It follows that
the analytic and dynamical components of $\Fat_{\phi}$ containing
$1$ are both the annulus $V$.  However, the D-component and
dynamical D-component are both the open disk $D(1,1)$.

In addition, for any $\alpha\in K$ with $|\alpha|=|c|$, write
$W_{\alpha}=D(\alpha,|c|)$.  There are infinitely many such disks
for which there are distinct integers $n>m\geq 0$
such that $\phi^m(W_{\alpha})=\phi^n(W_{\alpha})$.  For any
such $\alpha$, the analytic component, D-component, dynamical component,
and dynamical D-component all coincide and are equal to $W_{\alpha}$.
\end{example}

\section{Nontrivial reduction}
\label{sect:nontriv}

As is well known, the natural projection
$\ints_{\CK}\rightarrow \ints_{\CK}/\maxid_{\CK} = \khat$
induces a reduction map $\text{red}:\PCK\rightarrow\Pkhat$.
Given $\bar{a}\in\Pkhat$, the associated {\em residue class},
which we shall denote $W_{\bar{a}}\subset\PCK$, is the preimage
$$W_{\bar{a}} = \text{red}^{-1} ( \bar{a} ).$$
Any such class is either an open disk
$W_{\bar{a}}=D(a,1)$ with $a\in\CK$ and $|a|\leq 1$, or else
it is the disk at infinity, $W_{\bar{\infty}}=\PCK\setminus\Dbar(0,1)$.

Given a rational function $\phi\in\CK(z)$ and a residue class
$W_{\bar{a}}$, it will be useful to know whether or not
$\phi(W_{\bar{a}})$ is again a residue class.
To do so,
we recall the following definition of Rivera-Letelier \cite{Riv2},
which generalizes the notion of good reduction first
stated by Morton and Silverman \cite{MorSil1}.
\begin{defin}
\label{def:ntred}
Let $\phi\in \CK(z)$ be a nonconstant
rational function.  Write $\phi$ as $f/g$, with
$f,g\in\ints_{\CK}[z]$, such that at least one coefficient of
$f$ or $g$ has absolute value $1$.  Denote by
$\bar{f}$ and $\bar{g}$ the reductions of $f$ and $g$
in $\khat[z]$.  Let $\bar{h}=\gcd (\bar{f},\bar{g})\in \khat[z]$,
let $\bar{f}_0 = \bar{f}/\bar{h}$, and let
$\bar{g}_0 = \bar{g}/\bar{h}$.
We say that $\phi$ has {\em nontrivial
reduction} if $\bar{f}_0$ and $\bar{g}_0$ are not both
constant.  In that case, we define
$\bar{\phi} = \bar{f}_0/\bar{g}_0 \in \khat(z)$.
If $\deg \bar{\phi}=\deg\phi$,
we say $\phi$ has {\em good reduction}.
\end{defin}
If $\phi$ and $\psi$ have nontrivial reductions
$\bar{\phi}$ and $\bar{\psi}$,
then $\phi\circ\psi$ has nontrivial reduction
$\bar{\phi}\circ\bar{\psi}$.
Rivera-Letelier showed that the above definition of
good reduction is equivalent to Morton and Silverman's
original definition.
His analysis 
is summarized in the following two lemmas.  The
proofs, stated for the field $\Cp$, but which
apply to arbitrary $\CK$, appear
in \cite[Proposition 2.4]{Riv1}.
\begin{lemma}
\label{lem:hpfix}
Let $\phi\in \CK(z)$ be a rational function.
Then $\phi$ has nontrivial reduction if and only if
there are (not necessarily distinct) points
$\bar{a},\bar{b}\in\Pkhat$ such that
$\phi(W_{\bar{a}})=W_{\bar{b}}$.
\end{lemma}

\begin{lemma}
\label{lem:ntred}
Let $\phi\in \CK(z)$ be a rational function
of nontrivial reduction $\bar{\phi}\in\khat(z)$.  Then
there is a finite set $T\subset\Pkhat$ such that
$$\phi(W_{\bar{a}})= W_{\bar{\phi}(\bar{a})}
\quad\text{for all} \quad \bar{a}\in\Pkhat\setminus T,$$
and
$$\phi(W_{\bar{a}})= \PCK
\quad\text{for all} \quad \bar{a}\in  T.$$
Moreover, $\phi$ has good reduction if and only
if $T=\emptyset$.
\end{lemma}

Given $\phi\in\CK(z)$ of nontrivial reduction and
its set $T\subset\Pkhat$ from Lemma~\ref{lem:ntred},
we call classes $W_{\bar{a}}$ of elements $\bar{a}\in T$
the {\em bad classes}, and we call the remaining
classes the {\em good classes}.  The bad classes
are precisely those classes that contain {\em both} a zero and a
pole of $\phi$; that is, they are the classes $W_{\bar{a}}$
corresponding to linear factors $(z-\bar{a})$ of
$\bar{h}=\gcd (\bar{f},\bar{g})$ in Definition~\ref{def:ntred}.

%

The following lemma will be
needed to prove Theorem~\ref{thm:fixhpwd}.
We provide a sketch of the proof, using methods
similar to those used by Rivera-Letelier.
\begin{lemma}
\label{lem:fixclass}
Let $\phi\in K(z)$ be a rational function of
nontrivial reduction.  Let $\bar{a}\in\Pkhat$
be a point of ramification of $\bar{\phi}$ which
is also fixed by $\bar{\phi}$.
Let $0<r<1$, and let
$a\in\PCK$ be a point in the residue class $W_{\bar{a}}$.
If $\bar{a}\neq \infty$, let $U$ be the annulus
$D(a,1)\setminus\Dbar(a,r)$; if $\bar{a}=\infty$,
let $U$ be the image of $D(1/a,1)\setminus\Dbar(1/a,r)$
under the map $z\mapsto 1/z$.
Then the set
$$W_{\bar{a}} \setminus
\left( \bigcup_{n\geq 0} \phi^n(U) \right)$$
contains at most one point.
\end{lemma}

\begin{proof} (Sketch).
After a $\PGL(2,\ints_{\hat{K}})$-change of coordinates,
we may assume that $a=0$.  If $\bar{0}$ is a good class,
then the hypotheses imply that for $z\in D(0,1)$, $\phi(z)$
is given by a power series
$$\phi(z) = \sum_{i=0}^{\infty} c_i z^i$$
with all $|c_i|\leq 1$, with $|c_0|, |c_1| < 1$,
and with $|c_m|=1$ for some minimal $m\geq 2$.  (The conditions
on $c_1$ and $c_m$ come from the ramification hypothesis;
they imply that the reduction $\bar{\phi}$ looks like
$\bar{c}_m z^m + \bar{c}_{m+1}z^{m+1} + \ldots$.)
Solving $\phi(z)=z$, it follows easily
that $D(0,1)$ contains
a fixed point $b$; without loss, $b=0$, so that $c_0=0$.
Then for any $0<s<1$,
solving the power series equations $\phi(z)=x$
for $x\in D(0,1)\setminus\Dbar(0,s^m)$ shows that
$$ D(0,1)\setminus\Dbar(0,s^m) \subseteq
\phi\left( D(0,1)\setminus\Dbar(0,s) \right) .$$
Thus, for any nonzero $x\in D(0,1)$, there must be an
integer $n\geq 0$ such that $x\in\phi^n(U)$.
Hence, $D(0,1)\setminus\bigcup\phi^n(U)\subseteq\{0\}$.
(In dynamical language, $0$ is an attracting fixed point
with basin containing $D(0,1)$.)

If $\bar{0}$ is a bad class, then $D(0,1)$ contains
finitely many poles, so that for $z\in D(0,1)$, $\phi(z)$
may be written as
$$\phi(z) = \left(\sum_{i=0}^{\infty} c_i z^i\right)
+ \sum_{j=1}^M \frac{A_j}{(z-\alpha_j)^{e_j}}, $$
with the same conditions as before on $\{c_i\}$,
and with $|A_j|, |\alpha_j|<1$.
Again, we may change coordinates so that $c_0=0$,
although this time, $0$ itself may not be a fixed point.
Let $R=\max\{ |\alpha_j|\}<1$.
Then for any $s\in [R,1)$,
$$ D(0,1)\setminus\Dbar(0,s^m) \subseteq
\phi\left( D(0,1)\setminus\Dbar(0,s) \right) ,$$
as before.
Thus, $\phi^n(U)$ contains a pole for some $n\geq 0$;
further computations show that $\phi^{n+1}(U)=\PCK$.
\end{proof}

\section{Canonical heights}
\label{sect:height}

To prove our existence result (Theorem~\ref{thm:fixhpwd}),
we will need a few facts from the theory
of diophantine height functions.  We present
the required statements without proof; instead,
we refer the reader
to \cite[Chapters 2--4]{Lan} for more details.
The results of this section will be used only
in the technical proof of Lemma~\ref{lem:wandpt}.
The reader may therefore
wish to skip ahead to the application of
Lemma~\ref{lem:wandpt} in the proof
of Theorem~\ref{thm:fixhpwd}.

Let $k_0$ be either $\QQ$ or else the
field $L(T)$ of rational functions in one variable
defined over an arbitrary field $L$.
Let $k$ be a finite extension of $k_0$,
and let $\khat$ be an algebraic closure of $k$.


The standard {\em height function} $h_0:k_0\rightarrow\RR_{\geq 0}$
is given by
\begin{equation}
\label{eq:htdef1}
h_0\left( \frac{f}{g} \right) = \max\{ \deg f , \deg g\}
\end{equation}
if $k_0=L(T)$ and $f,g\in L[T]$ are relatively prime polynomials, or
\begin{equation}
\label{eq:htdef2}
h_0\left( \frac{m}{n} \right) = \max\{ \log |m| , \log |n| \}
\end{equation}
if $k_0=\QQ$ and $m,n\in\ZZ$ are relatively prime integers.
Considering $k_0$ as a subset of $\Pkhat$ in the natural way, the height
function $h_0$ extends to
$$h: \Pkhat \rightarrow \RR_{\geq 0}$$
with the property that for
any rational function $\bar{\phi}\in k(z)$
of degree $d\geq 1$,
there is a real constant $C=C_{\bar{\phi}}\geq 0$ such that
$$ \text{for all } x\in \Pkhat,
\quad \left| h\left(\bar{\phi}(x) \right) - d h(x)\right| \leq C .
$$

For a {\em fixed} function $\bar{\phi}(z)$ of degree $d\geq 2$,
Call and Silverman \cite{CS} introduced a related {\em canonical}
height function
$$\hat{h}=\hat{h}_{\bar{\phi}} : \Pkhat \rightarrow \RR_{\geq 0}$$
generalizing a construction of
N\'eron \cite{Ner} and Tate \cite{Tat}.
The key property of $\hat{h}$ is that there
is a real constant $C'=C'_{\bar{\phi}}\geq 0$ such that
for all $x\in\Pkhat$,
\begin{equation}
\label{eq:canht}
\hat{h}(\bar{\phi}(x) ) = d \hat{h}(x)
\quad \text{and} \quad
\left|\hat{h}(x) - h(x) \right| \leq C' ,
\end{equation}
where $h$ is the standard height function described above.
Note that by \eqref{eq:canht}, a preperiodic
point $x$ of $\bar{\phi}$ must have canonical height
$\hat{h}_{\bar{\phi}}(x)=0$.

The following lemma is not directly concerned
with heights, but it applies to fields $k$ of the type we
have been considering in this section.  It can be proven
using the fact such a field
contains a Dedekind ring of integers $\ints_k$ with infinitely
many prime ideals.
\begin{lemma}
\label{lem:cmap}
Let $k_0$ be either $\QQ$ or the function field $L(T)$ for some field $L$,
and let $k$ be an algebraic extension of $k_0$.
Let $c\in k^{*}$
such that $c^n\neq 1$ for all $n\geq 1$,
and define $\bar{\phi}(z) = cz$.
Then there exists an infinite sequence $\{x_i: i\in\ZZ\}$
of wandering points in $\Pk$ such that for any distinct $i,j\in\ZZ$,
$x_i$ and $x_j$ lie in different grand orbits of $\bar{\phi}$.
\end{lemma}

\section{Existence of wandering domains}
\label{sect:exist}

Our strategy for constructing wandering domains of $\phi(z)\in K(z)$
begins with finding wandering {\em points} in $\Pkhat$
of the reduction $\bar{\phi}(z)\in k(z)$.
The following
lemma shows that outside of trivial counterexamples,
such points always exist.
As mentioned in the previous section, the reader may
wish to skip the proof of the Lemma to see its
application in the proof of Theorem~\ref{thm:fixhpwd}.

\begin{lemma}
\label{lem:wandpt}
Let $k$ be a field, and let
$\bar{\phi}(z)\in k(z)$ be a nonconstant rational function.
Suppose that for every $n\geq 1$, $\bar{\phi}^n$ is not
the identity function.
Then the following five statements are equivalent:
\begin{list}{\rm \alph{bean}.}{\usecounter{bean}}
\item
$k$ is an algebraic extension of a finite field.
\item
Only finitely many wandering grand orbits
of $\bar{\phi}$ intersect $\Pk$.
\item
There are only finitely many wandering grand orbits
of $\bar{\phi}$ in $\Pkhat$.
\item
There are no wandering grand orbits
of $\bar{\phi}$ intersecting $\Pk$.
That is, all points in $\Pk$ are preperiodic under $\bar{\phi}$.
\item
There are no wandering grand orbits
of $\bar{\phi}$ in $\Pkhat$.
That is, all points in $\Pkhat$ are preperiodic under $\bar{\phi}$.
\end{list}
\end{lemma}

\begin{proof}
({\em i}).
Clearly (e) implies (d) implies (b),
and (e) implies (c) implies (b).

To show (a) implies (e),
suppose that $k$ is an algebraic extension of a finite field.
Then we may assume that $\khat \cong \Fpbar$, an
algebraic closure of the field $\Fp$ of $p$ elements,
for some prime number $p$.
Given $x\in\Pkhat$, there is some $r\geq 1$ such that
$x\in\PP^1(\FF_{p^r})$ and all of the (finitely many)
coefficients of $\bar{\phi}(z)$ also lie in $\FF_{p^r}$.
Since $\PP^1(\FF_{p^r})$ is a finite set which is
mapped into itself by $\bar{\phi}$, $x$ must
be preperiodic, proving the implication.


The remaining (and substantive) part of the proof is
to show that (b) implies (a).
Suppose that $k$ is not an algebraic extension
of a finite field; we must show that $\Pk$
intersects infinitely many wandering grand orbits of $\bar{\phi}$.

({\em ii}).
We will now reduce to the case that $k$ is
a finite extension either of $\QQ$ or
of the function field $L(T)$, for some field $L$.

Clearly, $k$ is a field extension of $L_0$,
where $L_0=\QQ$ if $\charact k=0$, or $L_0=\Fp$
if $\charact k=p>0$.  If $k/L_0$ is
an algebraic extension, then by hypothesis,
$k$ must be an algebraic extension of $\QQ$.

On the other hand, if
$k/L_0$ is a transcendental extension, then there
is a nonempty transcendence basis $B\subset k$ such
that $k$ is an algebraic extension of $L_0(B)$
(See, for example, \cite[Theorem~8.35]{Jac}.)
Pick $T\in B$, let $B'= B\setminus\{T\}$, and
let $L=L_0(B')$, so that $L(T) \cong L_0(B)$,
and $T$ is transcendental over $L$.
In that case, then, $k$ is an algebraic extension of the
function field $L(T)$.

We may now assume that $k$ is
a {\em finite} extension of either $\QQ$ or $L(T)$.
After all, the finitely many coefficients of $\bar{\phi}$ are
each algebraic over $\QQ$ or $L(T)$, so there is
a single finite extension that contains all of them.

Write $k_0=\QQ$ or $k_0=L(T)$ as appropriate,
so that $k$ is a finite extension of $k_0$.
Let $d=\deg\bar{\phi}$.  We consider two cases:
that $d=1$, or that $d\geq 2$.

({\em iii}).
If $d=1$, then by a change of coordinates, we may
assume that either $\bar{\phi}(z)=z+1$ or $\bar{\phi}(z)=cz$,
for some $c\in k^*$.
(If there are two distinct fixed points, move one to $0$
and one to $\infty$, to get $\bar{\phi}(z)=cz$.
If there is only one, move it to
$\infty$ and then scale to get $\bar{\phi}(z)=z+1$.)
If $\bar{\phi}(z)=z+1$ and $\charact k =p>0$,
then $\bar{\phi}^p(z)=z$, contradicting the hypotheses.
If $\bar{\phi}(z)=z+1$ and $\charact k = 0$,
then $\QQ\subset k$, so that there are clearly infinitely
many wandering grand orbits; for example, there is
one such orbit for each element of $\QQ\cap [0,1)$.
On the other hand, if $\bar{\phi}(z)=cz$, then by
hypothesis, $c^n\neq 1$ for all $n\geq 1$.  Therefore,
we have infinitely many wandering grand orbits
by Lemma~\ref{lem:cmap}.

({\em iv}).
For the remainder of the proof, suppose $d\geq 2$.
Define the canonical
height function $\hat{h}=\hat{h}_{\bar{\phi}}$
as in Section~\ref{sect:height}, and let $C'=C'_{\bar{\phi}}\geq 0$
be the corresponding constant in inequality~\eqref{eq:canht}.
Let $M=1 + 2C'$.

We claim that for any real number
$r\geq 0$, there exists $x\in\Pk$ satisfying
$\hat{h}(x)\in (r,r+M]$.
Indeed,
by the definition of the height function
$h_0$ in equations~\eqref{eq:htdef1} and~\eqref{eq:htdef2}, there is
some $x\in k_0$ such that
$h(x)\in (r+C',r+C'+1]$.
Because $|\hat{h}(x) - h(x)| \leq C'$,
it follows that $\hat{h}(x)\in (r, r+M]$,
proving the claim.

({\em v}).
Let $N\geq 1$ be any positive integer; we will
show that $\bar{\phi}$ has at least $N$ distinct
wandering grand orbits which intersect $\Pk$.

Let $I$ be the real interval $I=( MN, 2MN ]$.
By ({\em iv}), there are at least $N$ different
points $x\in\Pk$ such that $\hat{h}(x)\in I$.
Recall that the preperiodic points all 
have canonical height zero; so if
$\hat{h}(x)\in I$, then $x$ must be wandering.
Thus, it suffices to show that
if $x,y \in\Pk$ are two points
with $\hat{h}(x),\hat{h}(y)\in I$ but $\hat{h}(x)\neq\hat{h}(y)$,
then $x$ and $y$ must lie in different grand orbits.

Suppose not.  Then there exist points $x,y \in\Pk$
with $\hat{h}(x),\hat{h}(y)\in I$ but $\hat{h}(x)\neq\hat{h}(y)$,
and integers $n\geq m\geq 0$
such that $\bar{\phi}^m(x)=\bar{\phi}^n(y)$.
Thus, we have
$d^m\hat{h}(x) = d^n\hat{h}(y)$, and therefore
$\hat{h}(x) = d^{n-m}\hat{h}(y)$.
Since $\hat{h}(x)\neq\hat{h}(y)$, we must have $m<n$.
Hence,
$$ 2MN < 2 \hat{h}(y) \leq d^{n-m}\hat{h}(y) = \hat{h}(x)
\leq 2MN,$$
because $MN<\hat{h}(y),\hat{h}(x)\leq 2MN$.
This contradiction completes the proof.
\end{proof}

We are now prepared to state and prove our existence theorem,
which immediately implies part (a) of Theorem~A.
\begin{thm}
\label{thm:fixhpwd}
Let $K$ be a non-archimedean field
with residue field $k$, where $k$
is not an algebraic extension of
a finite field.
Let $\phi\in K(z)$ be a rational function
of nontrivial reduction $\bar{\phi}$,
and suppose that $\deg\bar{\phi}\geq 2$.
Then there is an infinite set
$\{\bar{b}_i : i\in\ZZ\}\subset\Pkhat$ such that
$\phi^n(W_{\bar{b}_i}) = W_{\bar{\phi}^n(\bar{b}_i)}$
for every $n\geq 0$ and $i\in\ZZ$, and such that all iterates
$\bar{\phi}^n(\bar{b}_i)$ are distinct.
Furthermore,
\begin{list}{\rm \alph{bean}.}{\usecounter{bean}}
\item
Each $W_{\bar{b}_i}$ is a wandering dynamical component
and a wandering dynamical D-component for $\phi$,
and each $W_{\bar{b}_i}$ lies in a different grand orbit
of such components.
\item 
If the Julia set $\calJ$ of $\phi$ intersects
at least two different residue classes
$W_{\bar{a}_1}$, $W_{\bar{a}_2}$, then
each $W_{\bar{b}_i}$ is also a wandering D-component,
and each $W_{\bar{b}_i}$ lies in a different grand orbit
of such components.
\item
If $\calJ$ has nonempty intersection
with infinitely many different residue classes,
then each $W_{\bar{b}_i}$ is a wandering analytic component for $\phi$,
and each $W_{\bar{b}_i}$ lies in a different grand orbit
of such components.
\end{list}
\end{thm}

\begin{proof}
({\em i}).
Let $\{\bar{c}_1,\ldots,\bar{c}_m\}\subset\Pkhat$
represent the finitely many
bad residue classes for $\bar{\phi}$.  We claim that there
is an infinite set $\{\bar{b}_i:i\in\ZZ\}\subset\Pkhat$
such that no $\bar{b}_i$ is preperiodic under $\bar{\phi}$,
such that $\bar{\phi}^n(\bar{b_i})$ avoids the $\bar{c}_j$'s,
and such that for any distinct $i,j\in\ZZ$, the grand
orbits of $\bar{b}_i$ and $\bar{b}_j$ under $\bar{\phi}$
are distinct.

To prove the claim, note that by Lemma~\ref{lem:wandpt},
there are points $\{\bar{b}'_i:i\in\ZZ\}$ in $\Pkhat$,
each with infinite forward orbit under $\bar{\phi}$,
such that no two lie in the same grand orbit.
For each $i\in\ZZ$,
let $N_i$ be the largest nonnegative integer $n$ such that
$\bar{\phi}^n(\bar{b}'_i)$ equals some $\bar{c}_j$,
or else $N_i=-1$ if no such $n$ exists.  Then
$\bar{b}_i=\bar{\phi}^{N_i+1}(\bar{b}'_i)$ for each $i$
satisfies the claim.

It follows immediately from Lemma~\ref{lem:ntred}
that for all $i\in\ZZ$ and all $n\geq 0$, 
$\phi^n(W_{\bar{b}_i})= W_{\bar{\phi}^n(\bar{b}_i)}$.
Thus, each $W_{\bar{b}_i}$ wanders and lies in the Fatou
set of $\phi$. 
Moreover, $W_{\bar{b}}$ is
a rational open disk, and therefore it must be contained in
a single component of the Fatou set, by any of the four
definitions of components.
Thus, it suffices only to show that each $W_{\bar{b}_i}$ is the full
Fatou component, for each of the four types.

({\em ii}).
Fix $\bar{b}=\bar{b}_i$ for some $i\in\ZZ$.
Let
$V_{dD}$ be the dynamical D-component containing $W_{\bar{b}}$,
$V_{dyn}$ the dynamical component,
$V_D$ the D-component, and
$V_{an}$ the analytic component.

If $V_{dyn} \supsetneq W_{\bar{b}}$, then
$V_{dyn}$ contains a connected open affinoid 
$U$ such that $U\supsetneq W_{\bar{b}}$.
Write $U=\PCK\setminus (D_1\cup\cdots\cup D_m)$
where $D_1,\ldots,D_m$ are disjoint closed disks.
Because $U$ {\em properly} contains a residue class,
each disk $D_i$ either is contained in $\PCK\setminus\Dbar(0,1)$
or has radius strictly less than $1$.
Thus, $U$ must 
contain all but finitely many
residue classes.  Define the finite (and possibly empty) sets
$$T_1=\left\{ \bar{a}\in\Pkhat : W_{\bar{a}}\not\subset U \right\}$$
and
$$T_2=\left\{ \bar{a}\in T_1 :
\bar{\phi}^{-n}(\bar{a}) \subseteq T_1
\text{ for all } n\geq 0 \right\}.$$
That is, $T_2$ is the set of all points $\bar{a}\in\Pkhat$
none of whose preimages $\bar{c}$ under any $\bar{\phi}^n$
have class $W_{\bar{c}}$ contained in $U$.
Because $T_2$ is finite and $\bar{\phi}^{-1}(T_2)\subset T_2$,
every element
of $T_2$ must be periodic under $\bar{\phi}$.
Thus, there is some integer $m\geq 1$
such that $\bar{\phi}^m$ fixes every element of $T_2$;
it follows that for every $\bar{a}\in T_2$,
$\bar{\phi}^{-m}(\bar{a}) = \{\bar{a}\}$.
But $\bar{\phi}^m$ has degree larger than $1$, and therefore every
element of $T_2$ is a fixed ramification point of $\bar{\phi}^m$.

Let $\tilde{U}=\bigcup_{n\geq 0}\phi^n(U)$.
For any $\bar{a}\not\in T_2$, there is some $\ell\geq 0$ such
that $\phi^{\ell}(U)\supset W_{\bar{a}}$;
hence $\tilde{U}\supset W_{\bar{a}}$.
On the other hand,
for $a\in T_2$, the intersection
$U\cap W_{\bar{a}}$ contains an annulus of the sort described
in Lemma~\ref{lem:fixclass}.  By that lemma, then,
$\tilde{U}$
must contain all but at most one point of $W_{\bar{a}}$.
Thus, $\tilde{U}$ contains all but finitely many points of $\PCK$,
which contradicts the definition of a dynamical component.
Therefore $V_{dyn}=W_{\bar{b}}$.

From the definitions,
we have $W_{\bar{b}}\subseteq V_{dD}\subseteq V_{dyn}$.
Thus, $V_{dD}=W_{\bar{b}}$ also.

({\em iii}).
Next, under the assumption that $\calJ$ intersects
at least two different residue classes,
suppose that $V_D\supsetneq W_{\bar{b}}$.
Then $V_D$ contains a disk $U\supsetneq W_{\bar{b}}$.  Such
a disk must contain all but one residue class,
and therefore $U$ must intersect the Julia set, which is
impossible.  Therefore $V_D=W_{\bar{b}}$.

Similarly, if $\calJ$ intersects
infinitely many classes $W_{\bar{a}}$,
and if $V_{an}\supsetneq W_{\bar{b}}$, then 
$V_{an}$ contains a connected affinoid $U\supsetneq W_{\bar{b}}$.
$U$ must contain all but finitely many residue classes,
which is impossible because then $U$ would intersect $\calJ$.
Hence $V_{an}= W_{\bar{b}}$.
\end{proof}

The following theorem shows that the condition that
the Julia set intersects infinitely many
different residue classes holds frequently.
\begin{thm}
\label{thm:julia}
Let $K$ be a non-archimedean field
with residue field $k$, let $p=\charact k \geq 0$,
let $\phi\in K(z)$ be a rational function
of nontrivial reduction $\bar{\phi}$,
and let $\Jul\subset\PCK$ be the Julia set of $\phi$.
Suppose either that $\bar{\phi}$ is separable
and of degree at least two, or
that there is a separable
map $\bar{\psi}\in k(z)$ of degree at
least two and an integer $r\geq 1$
such that $\bar{\phi}(z) = \bar{\psi}(z^{p^r})$.
If $\Jul$ intersects at least three different
residue classes of $\PCK$, then
$\Jul$ intersects infinitely many different
residue classes in $\PCK$.
\end{thm}

\begin{proof}
Because $\bar{\psi}$ is separable and of degree at least two,
then by the Riemann-Hurwitz formula
(see \cite[Corollary~2.4]{Har}), for example)
at most two points of $\Pkhat$
have only one preimage each under $\bar{\psi}$.
Given any $N\geq 3$, then, and any set $S_N\subset\Pkhat$
of $N$ distinct points, the number of points
in $\bar{\psi}^{-1}(S_N)$ must be strictly greater
than $N$.

Applying this fact inductively to $\bar{\phi}(z)=\bar{\psi}(z^{p^r})$,
we see that, given any
three distinct points $\bar{c}_1,\bar{c}_2,\bar{c}_3\in\Pkhat$,
there are infinitely many points $\bar{a}\in\Pkhat$ which
eventually map to some $\bar{c}_i$ under some $\bar{\phi}^n$.

Meanwhile, by Lemma~\ref{lem:ntred},
for any class $W_{\bar{a}}$, we have
$\phi(W_{\bar{a}})\supseteq W_{\bar{\phi}(\bar{a})}$.
Thus, if $\calJ$ intersects at least three residue
classes, it must intersect
infinitely many residue classes.
\end{proof}

To show that wandering domains coming from nontrivial
reduction actually exist, we present the following
example, which is just a generalization of
\cite[Example~2]{Ben5}.  Our example proves part~(b)
of Theorem~A.

\begin{example}
\label{ex:fixhpwd}
Let $K$ be a non-archimedean field
with residue field $k$ that is not
an algebraic extension of a finite field.
If $m\geq 2$ is an integer not divisible
by $\charact k$, let $\Psi_m(z)$ denote the $m$-th cyclotomic
polynomial.  For example, if $\charact k \neq 2$, we may
choose $m=2$ and hence $\Psi_m(z)=z+1$; if $\charact k = 2$,
we may choose $m=3$ and hence $\Psi_m(z)=z^2 + z +1$.  In
either case, $\bar{\Psi}_m$ has distinct roots in $\khat$, and if
$\zeta\in\CK$ is any root, then $\bar{\zeta}\neq 1$ but $\zeta^m=1$.

If $T\in K$ is any element satisfying $0< |T| < 1$,
define the rational function
$$\phi(z) = z^m + \frac{T}{\Psi_m(z)} =
	\frac{z^m\Psi_m(z) + T}{\Psi_m(z)}.$$
Then $\phi$ has nontrivial reduction $\bar{\phi}(z)=z^m\in\khat[z]$,
which is separable and of degree $m\geq 2$.  The only
bad residue classes are the roots of $\bar{\Psi}_m$ in $\khat$.
Hence, given $\bar{a}\in\Pkhat$ which is
not a root of $\bar{\Psi}_m$,
we have $\phi(W_{\bar{a}})=W_{\bar{a}^m}$.

Moreover, we claim that the Julia set of $\phi$ intersects
infinitely many distinct residue classes.  To show this,
let $\zeta\in \CK$ be a
root of $\Psi_m$.  First, we can easily check that $\phi$ has a fixed
point $\alpha\in\CK$ with $|\alpha-\zeta|=|T|$.  Indeed,
substituting $w=z-\zeta$ in the equation $\phi(z)=z$,
gives a polynomial in $\ints_{\CK}[w]$ with linear
coefficient $(1-\zeta)\Psi'_m(\zeta)$
(which has absolute value $1$) and constant term $T$
Second, we compute $|\phi'(\alpha)|=|T|^{-1} > 1$,
so that $\alpha$ is a repelling fixed point and hence lies in
the Julia set.  Furthermore, because $\bar{\phi}(z)=z^m$ is separable,
with no ramification points in $\Pkhat$ besides $0$ and $\infty$,
the set $\{\bar{\zeta}\}\cup \bar{\phi}^{-1}(\bar{\zeta})$
consists of at least three points.  Finally, the
corresponding residue classes each contain preimages
of $\alpha$, and hence they intersect the Julia set.
By Theorem~\ref{thm:julia}, our claim is valid.

It follows by Theorem~\ref{thm:fixhpwd} that
$\phi$ has infinitely
many grand orbits of wandering components
(of all four types).
More precisely, for any
$b\in K$ such that $\{\bar{b}^n\}_{n\in\ZZ}$ is an
infinite subset of $k$, the class $W_{\bar{b}}$ is a wandering
domain.  After all, no iterate $\bar{\phi}(\bar{b})$ is ever
one of the bad classes $\bar{\zeta}$, (otherwise, all future iterates
of $\bar{b}$ would be $\bar{1}$),
and those iterates are all distinct.
\end{example}

\smallskip

As another example, if the field $K$ satisfies the hypotheses
of Theorem~\ref{thm:fixhpwd}, then for $n\geq 2$,
the function $\phi(z)=z^n$
of Example~\ref{ex:comp1} has wandering dynamical
components and wandering dynamical D-components.
However, the unique analytic component
and D-component, namely $\PCK$, is not wandering.


As mentioned in the introduction,
sufficient conditions for residue
classes of $\PCK$ to be wandering domains
are more complicated if the
map has a nontrivial reduction of
degree one.  The remaining examples
of this section are of functions
of reduction degree one, all
defined over the field $K=\QQ((T))$,
whose residue field $\QQ$ is not
an algebraic extension of a finite
field.

\begin{example}
\label{ex:degonewd}
Let $\phi(z) = z + 1+ T/z \in \QQ((T))$.
Then $\phi$ has nontrivial
reduction $\bar{\phi}(z)=z+1$ of degree one.
The disk $D(0,1)$ contains the repelling fixed
point $-T$; 
it follows that $D(-m,1)$ intersects the
Julia set for every integer $m\geq 0$.
On the other hand,
the disk $U=D(1,1)$
satisfies $\phi^n(U)=D(n+1,1)$ for every
integer $n\geq 0$, so that $U$ lies in
the Fatou set.  Moreover, any strictly larger
affinoid containing $U$ must contain one
of the disks $D(-m,1)$ and hence must intersect
the Julia set.  So $U$ is a wandering
analytic component, wandering D-component,
wandering dynamical component,
and wandering dynamical D-component.
\end{example}

\begin{example}
\label{ex:degonenwd}
Let $\phi(z) = Tz^2 + z + 1 \in \QQ((T))$.
Then $\phi$ has nontrivial
reduction $\bar{\phi}(z)=z+1$ of degree one,
and as in the previous example,
the disk $U=D(1,1)$ lies in the Fatou
set and is wandering; in fact, the same
is true of every disk $D(b,1)$ for $|b|\leq 1$.
However, all these disks are contained
in the single disk $D(0,1/|T|)$, which
is fixed.  Thus, although the smaller
disks are wandering, none of them is large
enough to be a component of the Fatou set.

In fact, $\phi=h\circ\psi\circ h^{-1}$,
where $h(z) = Tz$ and $\psi(z)=z^2 + z + T$,
which is a map of good reduction, having
reduction $\bar{\psi}(z) = z^2 + z$
of degree two.  By Theorem~\ref{thm:fixhpwd},
$\psi$ does have wandering dynamical components
and wandering dynamical D-components.
Therefore $\phi$ also has such wandering
components, though they are {\em not} the residue classes $D(b,1)$
that we considered at first.
Moreover, the whole of $\PCK$ forms a single
D-component and a single analytic component;
hence, there are no wandering analytic or D-components.
\end{example}

\begin{example}
\label{ex:degonesemiwd}
Let $b=2$ and $c=T$ in Example~\ref{ex:comp5}, so that
$$\phi(z) = \frac{2z(z+T)(z+T^2)}{(z+2T)(z+T^3)(Tz+1)^2},$$
which has nontrivial
reduction $\bar{\phi}(z)=2z$ of degree one.
As in Example~\ref{ex:comp5}, let $V=D(0,|T|^{-1})\setminus\Dbar(0,|T|)$.
All of the residue
classes $D(a,1)$ (for $|a|=1$) are contained
in the Fatou set; in fact, every
such residue class $D(a,1)$
is a wandering D-component and a wandering dynamical D-component.
On the other hand, as we saw before, the affinoid $V$,
which contains all the disks $D(a,1)$,
is both a fixed analytic component and a fixed dynamical component.
\end{example}

\begin{example}
\label{ex:degonestat}
Let $b=-1$ and $c=T$ in Example~\ref{ex:comp5}, so that
$$\phi(z) = \frac{-z(z+T)(z+T^2)}{(z-T)(z+T^3)(Tz+1)^2},$$
which has nontrivial
reduction $\bar{\phi}(z)=-z$ of degree one.
Again, all of the residue
classes $D(a,1)$ (for $|a|=1$) are contained
in the Fatou set and are
both D-components and dynamical D-components.
This time, however, all those disks are
fixed by $\phi^2$, so none of them
is wandering.
As before, the open affinoid 
$V=D(0,|T|^{-1})\setminus\Dbar(0,|T|)$
contains all the disks $D(a,1)$
and is both a fixed analytic component and
fixed dynamical component.
\end{example}

\section{Residue characteristic zero}
\label{sect:czero}

We now prove Theorem~B.  The following theorem is a slightly
stronger result, showing that the desired conjugacy is defined
over a certain finite extension of $K$.

\begin{thm}
\label{thm:nowd}
Let $K$ be a discretely valued non-archimedean field with
residue field $k$ and residue characteristic
$\charact k =0$.
Let $\phi\in K(z)$ be a rational function,
and suppose that $U$ is a wandering analytic component,
wandering D-component, wandering dynamical D-component,
or wandering dynamical component
of $\phi$.  Let $L\subset\CK$ be any finite extension of $K$
such that $U$ contains a point of $\PP^1(L)$.
Then there
is a change of coordinates $g\in\PGL(2,L)$
and there are integers $M\geq 0$ and $N\geq 1$
such that
$\psi(z) =g\circ \phi^N\circ g^{-1}(z)$
has nontrivial reduction, $D(0,1)$ is
a wandering component (of the same type as $U$)
of $\psi$, and 
$g(\phi^M(U))\subset D(0,1)$.
\end{thm}

Note that a field $L$ satisfying the required properties
always exists.  Indeed, the algebraic closure $\hat{K}$
of $K$ is dense in $\PCK$, so that the open set
$U$ must contain some $a\in\hat{K}$.  Then
$L=K(a)$ is a finite extension of $K$.

\begin{proof}
We devote the bulk of the proof to the case that
$U$ is a wandering dynamical D-component.

({\em i}).
Let $a\in U\cap\PP^1(L)$.
Write $U_n=\Phi_{dD}^n(U)$ to simplify notation;
recall that $\Phi_{dD}^n(U)$ is the dynamical D-component containing
$\phi^n(U)$.  

We may assume without loss that $U_n\subset \Dbar(0,1)$
for every $n\geq 0$.  To do so, make a $\PGL(2,L)$-change
of coordinates to move $a$ to $\infty$ and $U$ to a set
containing $\PP^1\setminus\Dbar(0,1)$.  Because $U$
is wandering, it follows that $U_n\subset \Dbar(0,1)$
for every $n\geq 1$.  Finally, replace $U$ by $U_1$,
and we have the desired scenario.

Let $L'$ be a finite extension of $L$ such that
$\PP^1(L')$ contains all critical
points and all poles of $\phi$ in $\PCK$.
$L'$ is discretely valued, because it is
only a finite extension of the original field $K$.
Thus, there is a real number $0<\eps<1$ such that
$|(L')^*|=\{\eps^m : m\in\ZZ\}$.

Furthermore, there are only finitely many $n\geq 0$
such that $U_n$ contains a critical point,
and by Lemma~\ref{lem:dbound}, only finitely many $n\geq 0$
such that $\Phi_{dD}(U_n)\neq \phi(U_n)$.
Thus, by replacing $U$ by $U_{M'}$
for some $M'\geq 0$, we may
assume for all $n\geq 0$ that $U_n=\phi^n(U)$,
and that $U_n$ contains no critical points
and no poles.

Write $r_n=\diam (U_n)>0$, so that
$U_n=D(\phi^n(a),r_n)$, for each $n\geq 0$.
(Recall that the diameter and the radius of
a non-archimedean disk are the same; see
Section~\ref{ssec:nonarch}.)
By Lemmas~\ref{lem:isometry} and~\ref{lem:nocrit},
because each $U_n$ contains no critical points or poles,
there are integers $\ell_n\in\ZZ$ such that
$r_n=\eps^{\ell_n}r_0$.

({\em ii}).
We now claim that $r_0\in |(L')^*|$.
To prove the claim, suppose not.
Because $L'$ is discretely valued, there exists a real
number $s_0>r_0$ such that no $x\in L'$ satisfies $r_0\leq |x|< s_0$.
For every $n\geq 0$, let $s_n=r_n\cdot s_0/r_0$.
By the fact
that $r_n=\eps^{\ell_n}r_0$ and $|(L')^*|=\{\eps^m\}$,
it follows that no $x\in L'$ satisfies $r_n \leq |x| < s_n$.

Let $V_n=\phi^n(D(a,s_0))$, for all $n\geq 0$.
We will now show, by induction on $n$, that $V_n$ is an open
disk of radius (i.e., diameter)
$s_n$ that contains no critical points or poles.
For $n=0$, $V_0$
is an open disk of radius $s_0$ by definition,
and it contains no critical
points or poles because $V_0\cap L' = U_0\cap L'$ by
our choice of $s_0$.  Assuming the claim is true for $n\geq 0$,
then by Lemmas~\ref{lem:isometry} and~\ref{lem:nocrit},
$\diam(V_{n+1})/\diam(V_n) = \diam(U_{n+1})/\diam(U_n)$,
since $V_n$ contains no critical points or poles.
It follows that $V_{n+1}$ is a set of diameter $s_{n+1}$.
Thus, $V_{n+1}$ certainly omits at least two points of
$\PCK$; by Lemma~\ref{lem:diskmap}, then, it is an open disk.
Because no $x\in L'$ satisfies $r_{n+1} \leq |x| < s_{n+1}$,
we have $V_{n+1}\cap L' = U_{n+1}\cap L'$, and therefore
$V_{n+1}$ contains no critical points or poles, completing
the induction.

Since each $U_n$ is contained in $\Dbar(0,1)$,
then $s_n\leq s_0/r_0$ for every $n\geq 0$.
Therefore,
$$\phi^n(V_0) = V_n \subseteq \Dbar(0,s_0/r_0)$$
for all $n\geq 0$.
Because $U=U_0\subsetneq V_0$, we have contradicted
the assumption that $U$ is a dynamical D-component.

Thus, $r_0\in |(L')^*|$, as claimed.
It follows that $r_n\in |(L')^*|$ for all $n\geq 0$.

({\em iii}).
For all $n\geq 0$, let $\bar{U}_n=\Dbar(\phi^n(a),r_n)$, so that
$U_n \subsetneq \bar{U}_n\subset \Dbar(0,1)$.

We claim that for infinitely many
$n\geq 0$, $\bar{U}_n$ contains a pole or a critical point of $\phi$.
To prove the claim, suppose only finitely many of the $\bar{U}_n$
contained poles or critical points,
and replace $U$ by $U_{M'}$ (for some appropriate $M'\geq 0$)
so that no $\bar{U}_n$ contains a pole or critical point.
Because $|(L')^*|=\{\eps^m\}$ is discrete and $r_n\in |(L')^*|$,
the larger disk $D(\phi^n(a),r_n/\eps)$
also contains no poles or critical points for any $n\geq 0$.

For all $n\geq 0$, define
$V'_n = \phi^n(D(a,r_0/\eps))$.
By an induction argument similar to that in part ({\em ii}) above,
$V'_n = D(\phi^n(a),r_n/\eps)$.
Because $V'_0$ is an open disk that properly contains $U$,
and because $\phi^n(V'_0)\subseteq D(0,1/\eps)$,
we have contradicted
the assumption that $U$ is a dynamical
D-component, thus proving the claim.

({\em iv}).
Next, we claim that either $\bar{U}_n$ contains a pole
for infinitely many $n\geq 0$, or else 
there exist $M\geq 0$ and $N\geq 1$ such that
$\bar{U}_M = \bar{U}_{M+N}$.

If only finitely many of the $\bar{U}_n$ contain poles,
then by ({\em iii}), infinitely many of them contain critical
points.  As there are only finitely many critical points,
there must be integers $M\geq 0$ and $N\geq 1$ such that
$\bar{U}_M \cap \bar{U}_{M+N}$ is nonempty,
and such that for all $n\geq M$, $\bar{U}_n$ contains no poles.
Replacing $U$ by $U_M$ and $\phi$
by $\phi^N$, we may assume that $M=0$ and $N=1$.
By Lemma~\ref{lem:opcl} and the fact that $\bar{U}_n$
contains no poles, $\phi(\bar{U}_n)=\bar{U}_{n+1}$
for all $n\geq 0$.
Because $\bar{U}_0$ and $\bar{U}_1$ are disks
in $\CK$, either
$\bar{U}_0 \supsetneq \bar{U}_1$ or
$\bar{U}_0 \subseteq \bar{U}_1$.

If $\bar{U}_0 \supsetneq \bar{U}_1$,
then because $|(L')^*|=\{\eps^m\}$, we must have $r_1\leq \eps r_0 < r_0$.
Let $V'' = D(\phi(a),r_0)$.
Since $V''\subset \bar{U}_0$, we have
$\phi(V'')\subset \phi(\bar{U}_0)=\bar{U}_1\subset \bar{U}_0$;
by induction, we get
$\phi^n(V'')\subset \bar{U}_0$ for all $n\geq 0$.
However, $V''$ is an open disk that properly contains
$U_1$, contradicting the supposition that $U_1=\Phi_{dD}(U)$ is
a dynamical D-component.

If $\bar{U}_0 \subseteq \bar{U}_1$, then by the fact
that $\phi(\bar{U}_n)=\bar{U}_{n+1}$
for every $n\geq 0$, we have
$$\bar{U}_0 \subseteq \bar{U}_1 \subseteq \bar{U}_2 \subseteq \cdots .$$
If all the inclusions are proper, then
$$r_0 < r_1 < r_2 <\cdots .$$
Because $|(L')^*|=\{\eps^m\}$, we must have $r_n > 1$
for some $n\geq 0$, contradicting the assumption that
every $U_n$ is contained in $\Dbar(0,1)$.
Thus, for some $n\geq 0$, we have $\bar{U}_n=\bar{U}_{n+1}$,
and the claim is proven.
(In fact, we would have $\bar{U}_0=\bar{U}_1$,
but we do not need that result here.)

({\em v}).
Consider the case that $\bar{U}_n$ contains a pole
for infinitely many $n\geq 0$.
Since there are only finitely many poles,
there must be some pole $y\in\PP^1(L')$
and an infinite set $I$ of nonnegative integers
such that $y\in\bar{U}_n = \Dbar(\phi^n(a),r_n)$
for all $n\in I$.  Pick $s>0$ so that
$\phi(D(y,s))\subset \PCK\setminus \Dbar(0,1)$.  By
our initial assumptions, no $U_n$ can intersect $D(y,s)$,
or else $U_{n+1}$ would not be contained in $\Dbar(0,1)$.
Thus, $s<r_n\leq 1$  for all $n\in I$.

However, we also know that $r_n=\eps^{\ell_n}r_0\in |(L')^*|$
for all $n\geq 0$.  As $n$ ranges over $I$, then,
there are only finitely many possible values that $r_n$
can attain.  At least one must be attained
infinitely often.
In particular, there are integers
$M\geq 0$ and $N\geq 1$ such that $M, M+N \in I$
and $r_M=r_{M+N}$.  Since $y$ lies in both $U_M$
and $U_{M+N}$, we have $\bar{U}_M=\bar{U}_{M+N}$.

({\em vi}).
By ({\em iv}) and ({\em v}), then, there exist integers
$M\geq 0$ and $N\geq 1$ such that $\bar{U}_M=\bar{U}_{M+N}$.
Thus, $r_M=r_{M+N}$,
and $|\phi^M(a) -\phi^{M+N}(a)|\leq r_M$.  Because
$U_M\cap U_{M+N}=\emptyset$, we must in fact
have $|\phi^M(a) -\phi^{M+N}(a)|= r_M$.  Therefore
$r_M\in |L^*|$, since $a\in\PP^1(L)$ and $\phi\in K(z)\subseteq L(z)$.

Let $g\in\PGL(2,L)$ be the unique
linear fractional transformation satisfying
$g(\infty)=\infty$, $g(\phi^M(a))=0$, and $g(\phi^{M+N}(a))=1$.
Thus, $g(U_M)=D(0,1)$ and $g(U_{M+N})=D(1,1)$.
Let $\psi= g\circ \phi^N \circ g^{-1}$.
By Lemma~\ref{lem:hpfix}, $\psi$ has nontrivial reduction,
and the remaining conclusions of the theorem follow as well,
at least for the case of dynamical D-components.

({\em vii}).
Finally, suppose that $U$ is a wandering D-component,
wandering analytic component, or wandering dynamical
component containing a point $a\in\PP^1(L)$.
Let $U'$ be the dynamical D-component containing $a$.
Then $U'\subseteq U$, and therefore $U'$ is wandering.

For any integer $n\geq 0$,
define $U_n=\Phi^n(U)$ 
(where $\Phi$ is $\Phi_D$, $\Phi_{an}$, or $\Phi_{dyn}$,
as appropriate), and define $U'_n=\Phi_{dD}(U')$.
Choose $g$, $M$, and $N$ for $\phi$ as
in the theorem applied to $U'$.
It suffices to show that $U_M= U'_M$.

Suppose not; then $U'_M\subsetneq U_M$,
and therefore $g(U'_M)\subsetneq g(U_M)$.
Thus, $g(U_M)$ contains an affinoid strictly containing
the residue class $g(U'_M)$; hence, $g(U_M)$ contains
all but finitely many of the residue classes $D(b,1)$.
However, for every $n\geq 0$, $g(U'_{M+nN})$ is a residue
class.
In particular, $g(U_M)$ contains $g(U'_{M+nN})$
for some $n\geq 1$.  Thus, $U_M\cap U_{M+nN}$ is
nonempty, contradicting the wandering assumption
and proving the theorem.
\end{proof}

\appendix
\section{Rational functions and non-archimedean analysis}

\subsection{General properties of rational functions}
We recall some basic facts about rational functions $\phi\in L(z)$,
for an abstract field $L$ with algebraic closure $\hat{L}$.
A point $x\in L$ is called a {\em pole} of $\phi$ if
$\phi(x)=\infty$.
We may define the derivative $\phi'(z)$ away from
the poles by the usual formal
differentiation rules; if $L$ has a metric structure, then
the formal definition of $\phi'$ agrees with the difference
quotient definition of $\phi'$.

If $x\in\PP^1(\hat{L})$ maps to $\phi(x)$ with multiplicity greater than one
(i.e., if $\phi'(x)=0$), we say $x$ is a {\em critical point}
or {\em ramification point} of $\phi$.  After a coordinate
change in the domain and range, we may assume that $x=\phi(x)=0$,
and we may expand $\phi$ locally about $0$ as a power series
$$\phi(z) = \sum_{n=1}^{\infty} c_n z^n.$$
We say that {\em $x$ maps to $\phi(x)$ with multiplicity $m$} if
$m$ is the smallest integer such that $c_m\neq 0$.  Note that
if $\charact L =p>0$, the multiplicity might {\em not} be
the same as the number of the first nonzero derivative at $x$.
For example, if $\phi(z)=z^p$ where $\charact L=p$,
then $\phi'(z)=0$, but every
point $x$ maps to its image with multiplicity $p$, not
infinite multiplicity.

If $\phi'(z)$ is not identically zero, we
say $\phi$ is {\em separable}.  If $\charact L=0$, then
all nonconstant rational functions are separable.  If
$\charact L=p>0$, then $\phi\in L(z)$ is separable
if and only if $\phi$ cannot be written as
$\phi(z)=\psi(z^p)$ for any $\psi\in L(z)$.

A function $\phi\in L(z)$ may be written
as $\phi=f/g$, where $f,g\in L[z]$ are relatively prime polynomials.
The {\em degree} $\deg\phi$ is defined to be
$$\deg\phi = \max\{\deg f, \deg g\}.$$
Any point $y\in\PP^1(\hat{L})$ has exactly $\deg\phi$ preimages
in $\phi^{-1}(y)$, counting multiplicity.
If $\phi$ is separable of degree $d$,
then $\phi$ has exactly $2d-2$ critical points in $\PP^1(\hat{L})$,
counting multiplicity.  (Here, the multiplicity of a critical
point $x$ is the multiplicity of $x$ as a root of the equation
$\phi'(z)=0$.  Usually, this multiplicity is $e_x-1$, where $x$
maps to $\phi(x)$ with multiplicity $e_x$.  However, if
$\charact L=p>0$, and if $p\mid e_x$, then the multiplicity
of $x$ as a critical point will be strictly greater than $e_x -1$.
See \cite[IV.2]{Har} for more details.)

\subsection{Non-archimedean analysis}
\label{ssec:nonarch}
Given $a\in\CK$ and $r>0$, we denote by $D(a,r)$ and by
$\Dbar(a,r)$ the open disk and the closed disk, respectively,
of radius $r$ centered at $a$.  (We will follow the convention
that all disks have {\em positive} radius by definition, so
that singleton sets and the empty set are {\em not} considered
to be disks.)  By the non-archimedean triangle
inequality, any point of such a disk may be considered a center,
and if $U_1, U_2\subset\CK$ are two overlapping disks, then either
$U_1\subseteq U_2$ or $U_2\subseteq U_1$.  Moreover, if
$U\subset\CK$ is an open or closed disk of radius $r$,
then $r$ is also the diameter of $U$; that is,
$$r=\diam(U) = \sup \{|x-y| : x,y\in U \}.$$

The set $|K^*|=\{|x|: x\in K\setminus\{0\} \}\subset\RR_{>0}$
may be a discrete subset of $\RR_{>0}$; if so, we say
that $K$ is {\em discretely valued}.  In that case,
there is a real number $0<\eps < 1$ such that
$|K^*|=\{\eps^m : m\in \ZZ\}$.

Meanwhile, the set $|\CK^*|=\{|x|: x\in\CK\setminus\{0\} \}$ must be
dense in $\RR_{>0}$, but it need not contain
all positive real numbers.  For example,
$|\Cp^*|=\{p^q: q\in\QQ\}$.  Therefore, we say that
a disk $U$ is {\em rational} if $\diam(U)\in |\CK^*|$,
and $U$ is {\em irrational} otherwise.
If $a\in\CK$ and $r\in|\CK^*|$, then $D(a,r)\subsetneq \Dbar(a,r)$;
but if $r\in (\RR_{>0} \setminus |\CK^*|)$, then
$D(a,r)=\Dbar(a,r)$.  Thus, every disk is exactly one
of the following three types:
rational open, rational closed, or irrational.  The distinctions
between the three indicate metric properties, but not
topological properties; 
{\em all} disks are both open and closed as
{\em topological sets}.

More generally, a set $U\subset\PCK$ is a
{\em rational open disk} if either $U\subset\CK$ is a rational
open disk or $\PCK\setminus U$ is a rational closed disk.
Similarly, $U\subset\PCK$ is a
{\em rational closed disk} if either $U\subset\CK$ is a rational
closed disk or $(\PCK\setminus U)\subset\CK$
is a rational open disk;
and $U\subset\PCK$ is an
{\em irrational disk} if either $U\subset\CK$ is an irrational
disk or $(\PCK\setminus U)\subset\CK$
is an irrational disk.  There is a natural spherical metric
on $\PCK$ (see, for example, \cite{Ben4,Ben7,MorSil2}), but not all
the disks we have just defined in $\PCK$ are disks with
respect to the spherical metric.

If $U_1,U_2\subset\PCK$ are disks such that $U_1\cap U_2\neq\emptyset$
and $U_1\cup U_2\neq \PCK$, then either $U_1\subseteq U_2$
or $U_2\subseteq U_1$.  In particular, both $U_1\cap U_2$
and $U_1\cup U_2$ are also disks; and if
$U_1$ and $U_2$ are both rational closed (respectively,
rational open, irrational), then so are
$U_1\cap U_2$ and $U_1\cup U_2$.

The group $\PGL(2,\CK)$ acts on $\PCK$ by linear fractional
transformations.  Any $g\in\PGL(2,\CK)$ maps
rational open disks to rational open disks,
rational closed disks to rational closed disks, and
irrational disks to irrational disks.

The following four lemmas concern the action of
non-archimedean rational functions on disks.  In fact,
all four apply more generally to power series on disks,
though we do not need to define the necessary terminology
of rigid analyticity to state the lemmas.
We omit the proofs,
which are easy applications of the Weierstrass Preparation Theorem,
Newton polygons,
and other fundamentals of non-archimedean analysis.  Some
proofs may be found in \cite{Ben7}; see any of
\cite[Chapter 5]{BGR}, \cite[Chapter II]{Esc},
\cite[Chapter IV]{Kob}, or \cite[Chapter 6]{Rob}
for the theory surrounding such results.

\begin{lemma}
\label{lem:diskmap}
Let $U\subset\PCK$ be a disk, and let $\phi\in\CK(z)$ be a
rational function.
Suppose that $\PCK\setminus \phi(U)$ contains at
least two points.  Then $\phi(U)$ is a disk of
the same type (rational closed, rational open,
or irrational) as $U$.
\end{lemma}

\begin{lemma}
\label{lem:opcl}
Let $a,b\in\CK$,
let $r,s>0$,
and let $\phi\in\CK(z)$ be a
rational function with no poles in $\Dbar(a,r)$,
such that $\phi(D(a,r))=D(b,s)$.
Then $\phi(\Dbar(a,r))=\Dbar(b,s)$.
\end{lemma}

\begin{lemma}
\label{lem:isometry}
Let $U\subset\CK$ be a disk, let $a\in U$,
and let $\phi\in\CK(z)$ be a rational function
with no poles in $U$.
Then the following two statements are equivalent.
\begin{list}{\rm \alph{bean}.}{\usecounter{bean}}
\item
$\phi$ is one-to-one on $U$.
\item
For all $x,y\in U$, $|\phi(x)-\phi(y)| = |\phi'(a)|\cdot |x-y|$.
\end{list}
\end{lemma}

\begin{lemma}
\label{lem:nocrit}
Let $K$ be a non-archimedean field with residue
field $k$, and suppose that $\charact k=0$.
Let $U\subset\CK$ be a disk,
and let $\phi\in\CK(z)$ be a
rational function.
Then the following two statements are equivalent.
\begin{list}{\rm \alph{bean}.}{\usecounter{bean}}
\item
$\phi$ is one-to-one on $U$.
\item
$\phi$ has no critical points in $U$.
\end{list}
\end{lemma}

Lemma~\ref{lem:nocrit} is needed only in parts~(ii) and~(iii)
of the proof of Theorem~\ref{thm:nowd}, and it is the only
use in that proof of the hypothesis that $\charact k =0$.
(The lemma is also quoted in part~(i) of the same proof,
but its use there can be avoided if desired.)

\subsection{Rigid Analysis}
\label{ssec:rigid}
We will need some basic facts and definitions from the
non-archimedean theory of rigid analysis.  We refer
the reader to \cite[Part C]{BGR} or \cite{FvP} for detailed background,
or to \cite{Goss} for a broader (but still technical)
overview of the subject; however, the discussion that
follows is mostly self-contained.

A {\em connected affinoid} is a set
$W\subset\PCK$ of the form
$$W=\PCK \setminus \left( U_1 \cup U_2 \cup \cdots \cup U_N \right),$$
where $N\geq 0$, and where the $\{U_i\}$ are pairwise disjoint disks.
If each $U_i$ is rational open, we say $W$
is a {\em connected rational closed affinoid};
if each $U_i$ is rational closed, we say $W$
is a {\em connected rational open affinoid};
and 
if each $U_i$ is irrational, we say $W$
is a {\em connected irrational affinoid}.

If $W_1$ and $W_2$ are connected affinoids, and if
$W_1\cap W_2\neq \emptyset$, then $W_1\cap W_2$
and $W_1\cup W_2$ are also connected affinoids.
In that case, if
$W_1$ and $W_2$ are both rational closed (respectively,
rational open, irrational), then so are
$W_1\cap W_2$ and $W_1\cup W_2$.

In general, an {\em affinoid} is a finite union of connected
affinoids.  However, we will not need that notion in this paper.
Note that our definition allows the full set $\PCK$ and
the empty set $\emptyset$ to be considered
connected affinoids, while traditional rigid analysis does not.
Also note that we consider $\PCK$ to be a connected affinoid
of all three types.  Every other connected affinoid is at most
one of the three types; or, it may be none of them, if, for example,
$U_1$ is a rational open disk and 
$U_2$ is a rational closed disk.

Intuitively, connected affinoids are supposed to behave like
connected sets, even though {\em topologically}, all subsets of $\PCK$
are totally disconnected.  In particular, it is well known
(as can be shown using standard rigid analysis techniques)
that if $\phi\in\CK(z)$ is a rational function of degree $d$,
and if $W$ is a connected affinoid, then:
\begin{itemize}
\item $\phi(W)$ is also a connected affinoid.
Moreover, if $W$ is rational closed
(respectively, rational open, irrational), then so is $\phi(W)$.
\item $\phi^{-1}(W)$ is a disjoint union of
connected affinoids $V_1,\ldots, V_N$, with $1\leq N\leq d$.
For every $i=1,\ldots, N$,
$\phi$ maps $V_i$ onto $W$.  Moreover, if $W$ is rational closed
(respectively, rational open, irrational), then so are
$V_1,\ldots, V_N$.
\end{itemize}

The following lemma shows that for any given rational function
$\phi$, most disks $U\subset\PCK$ have preimage $\phi^{-1}(U)$
consisting simply of a finite union of disks.  It appeared
as \cite[Lemma 3.1.4]{Ben1}, but we include a partial proof here
for the convenience of the reader.

\begin{lemma}
\label{lem:dbound}
Let $U_1,\ldots, U_n\subset\PCK$ be
disjoint disks, and let $\phi\in K(z)$ be a rational
function of degree $d\geq 1$.
Suppose that for each $i=1,\ldots, n$, the
inverse image $\phi^{-1}(U_i)$ is not a finite
union of disks.  Then $n\leq d-1$.
\end{lemma}

\begin{proof}
(Sketch).
If $U_i$ is an open disk, then it can be
written as a nested union $\bigcup_{m\geq 1} V_m$ of rational closed
disks $V_m$, with $V_m\subset V_{m+1}$.  If each $\phi^{-1}(V_m)$
is a union of at most $d$ disks, then the same is true of
$\phi^{-1}(U_i)$.  Thus, we may assume that each $U_i$ is
a rational closed disk.

Let $W=\PCK\setminus(U_1\cup \cdots \cup U_n)$.  Then $W$
is a rational open connected affinoid.  By the discussion above,
the inverse image $\phi^{-1}(W)$ is a disjoint union of at
most $d$ rational open connected affinoids.  Thus, as we leave
to the reader to verify, the complement $\PCK\setminus\phi^{-1}(W)$
is a union of some rational closed disks and at most $d-1$
connected affinoids which are not disks.  (Note, for example,
that if $V$ is a closed affinoid that is neither a disk
nor all of $\PCK$, then the complement of $V$ consists of
at least two connected components.  Thus, if
$\PCK \setminus \phi^{-1}(W)$ consisted of $d$ or more
non-disk connected components, then $\phi^{-1}(W)$ would
consist of at least $d+1$ connected components.)
However, the complement
of $\phi^{-1}(W)$ is precisely the disjoint union
$\bigcup_{i=1}^n \phi^{-1}(U_i)$.  It follows that
$n\leq d-1$.
\end{proof}

\bibliographystyle{plain}

\begin{thebibliography}{99}

\bibitem{Bea}
A.~Beardon,
\newblock{{\em Iteration of Rational Functions},}
\newblock{Springer-Verlag, New York, 1991.}

\bibitem{Ben1}
R.~Benedetto,
\newblock{{\em Fatou components in $p$-adic dynamics},}
\newblock{Ph.D. thesis, Brown University, 1998.}

\bibitem{Ben2}
R.~Benedetto,
\newblock{$p$-adic dynamics and Sullivan's No Wandering Domains Theorem,}
\newblock{{\em Compositio Math.\/} {\bf 122} (2000), 281--298.}

\bibitem{Ben3}
R.~Benedetto,
\newblock{Hyperbolic maps in $p$-adic dynamics,}
\newblock{{\em Ergodic Theory Dynam. Systems\/} {\bf 21} (2001), 1--11.}

\bibitem{Ben4}
R.~Benedetto,
\newblock{Reduction, dynamics, and Julia sets of rational functions,}
\newblock{{\em J. Number Theory} {\bf 86} (2001), 175--195.}

\bibitem{Ben5}
R.~Benedetto,
\newblock{Components and periodic points in non-archimedean dynamics,}
\newblock{{\em Proc. London Math. Soc.} (3) {\bf 84} (2002), 231--256.}

\bibitem{Ben6}
R.~Benedetto,
\newblock{Examples of wandering domains in $p$-adic polynomial dynamics,}
\newblock{{\em C.~R.~Acad. Sci. Paris, Ser. I} {\bf 335} (2002), 615--620.}

\bibitem{Ben7}
R.~Benedetto,
\newblock{Non-archimedean holomorphic maps and the Ahlfors
Islands Theorem,}
\newblock{{\em Amer. J. Math.}, {\bf 125} (2003), 581--622.}

\bibitem{Ben8}
R.~Benedetto,
\newblock{Wandering domains in non-archimedean polynomial dynamics,}
\newblock{preprint, 2004.  Available online at
  {\tt http://arxiv.org/abs/math.NT/0312029} }

\bibitem{Ber}
V.~Berkovich,
\newblock{{\em Spectral theory and analytic geometry over
    non-Archimedean fields,\/}}
\newblock{Amer.\ Math.\ Soc., Providence, 1990.}

\bibitem{Bez}
J.-P.~B\'{e}zivin,
\newblock{Sur les points p\'{e}riodiques des applications
	rationnelles en analyse ultram\'{e}trique,}
\newblock{{\em Acta Arith.} {\bf 100} (2001), 63--74.}

\bibitem{BGR}
S.~Bosch, U.~G\"{u}ntzer, and R.~Remmert,
\newblock{{\em Non-Archimedean Analysis: A Systematic Approach to
        Rigid Analytic Geometry,\/}}
\newblock{Springer-Verlag Berlin, 1984.}

\bibitem{CS}
G.~Call and J.~Silverman,
\newblock{Canonical heights on varieties with morphisms,}
\newblock{{\em Compositio Math.\/} {\bf 89} (1993), 163--205.}

\bibitem{CG}
L.~Carleson and T.~Gamelin,
\newblock{ \em Complex Dynamics,\/}
\newblock{Springer-Verlag, New York, 1991.}

\bibitem{Esc}
A.~Escassut,
\newblock{{\em Analytic Elements in $p$-adic Analysis,\/}}
\newblock{World Scientific, Singapore, 1995.}

\bibitem{FvP}
J.~Fresnel and M.~van~der~Put,
\newblock{{\em G\'{e}om\'{e}trie Analytique Rigide et Applications,\/}}
\newblock{PM 18, Birkh\"{a}user, Boston, 1981.}

\bibitem{Goss}
D.~Goss,
\newblock{A short introduction to rigid analytic spaces, in}
\newblock{{\em The arithmetic of function fields}, Columbus, OH, 1991,
        131--141.}

\bibitem{Gou}
F.~Gouv\^{e}a,
\newblock{{\em $p$-adic Numbers. An Introduction,} 2nd ed.,}
\newblock{Springer-Verlag, Berlin, 1997.}

\bibitem{Har}
R.~Hartshorne,
\newblock{{\em Algebraic Geometry},}
\newblock{Springer-Verlag, New York, 1977.}

\bibitem{HY}
M.~Herman and J.-C.~Yoccoz,
\newblock{Generalizations of some theorems of small divisors to
        non-Archimedean fields, in}
\newblock{{\em Geometric Dynamics (Rio de Janairo, 1981),\/}
        LNM 1007, Springer-Verlag, Berlin-New York, 1983,
        pp.~408--447.}

\bibitem{Hs2}
L.-C.~Hsia,
\newblock{Closure of periodic points over a nonarchimedean field,}
\newblock{{\em J. London Math. Soc.} (2) {\bf 62} (2000), 685--700.}

\bibitem{Jac}
N.~Jacobson,
\newblock{{\em Basic Algebra II},}
\newblock{2nd ed., W.H.~Freeman, New York, 1989.}

\bibitem{Kiwi}
J.~Kiwi,
\newblock{Puiseux series polynomial dynamics and
  iteration of complex cubic polynomials,}
\newblock{preprint, 2004.  Available online at
  {\tt http://arxiv.org/abs/math.DS/0409397}}

\bibitem{Kob}
N.~Koblitz,
\newblock{{\em $p$-adic Numbers, $p$-adic Analysis, and Zeta-Functions},}
\newblock{2nd ed., Springer-Verlag, New York, 1984.}

\bibitem{Lan}
S.~Lang,
\newblock{{\em Fundamentals of Diophantine Geometry},}
\newblock{Springer-Verlag, New York, 1983.}

\bibitem{Mil}
J.~Milnor,
\newblock{{\em Dynamics in One Complex Variable: Introductory Lectures,}
	2nd ed.,}
\newblock{Vieweg, Braunschweig, 2000.}

\bibitem{MorSil1}
P.~Morton and J.~Silverman,
\newblock{Rational periodic points of rational functions,}
\newblock{{\em Inter. Math. Res. Notices\/} {\bf 2} (1994), 97--110.}

\bibitem{MorSil2}
P.~Morton and J.~Silverman,
\newblock{Periodic points, multiplicities, and dynamical units,}
\newblock{{\em J. Reine Angew. Math.\/} {\bf 461} (1995), 81--122.}

\bibitem{Ner}
A.~N\'eron,
\newblock{Quasi-fonctions et hauteurs sur les vari\'et\'es ab\'eliennes,}
\newblock{{\em Ann. of Math.} (2) {\bf 82} (1965), 249--331.}

\bibitem{Riv1}
J.~Rivera-Letelier,
\newblock{{\em Dynamique des fonctions rationnelles sur des corps locaux},}
\newblock{Ph.D. thesis, Universit\'{e} de Paris-Sud, Orsay, 2000.}

\bibitem{Riv2}
J.~Rivera-Letelier,
\newblock{Espace hyperbolique $p$-adique et dynamique de
	fonctions rationnelles,}
\newblock{{\em Compositio Math.} {\bf 138} (2003), 199--231.}

\bibitem{Rob}
A.~Robert,
\newblock{{\em A Course in $p$-adic Analysis},}
\newblock{Springer-Verlag, New York, 2000.}

\bibitem{RuBa}
R.~Rumely and M.~Baker,
\newblock{Analysis and dynamics on the Berkovich projective line,}
\newblock{preprint, 2004.  Available online at
  {\tt http://arxiv.org/math.NT/0407433}}

\bibitem{Ser2}
J.-P.~Serre,
\newblock{{\em Local Fields}, translated by M.~Greenberg,}
\newblock{Springer-Verlag, New York-Berlin, 1979.}
\newblock{(Translation of {\em Corps Locaux}, 2nd ed.,
Hermann, Paris, 1968.)}

\bibitem{SW}
N.~Smart and C.~Woodcock, 
\newblock{$p$-adic chaos and random number generation,}
\newblock{{\em Experiment. Math.\/} {\bf 7} (1998), no.~4, 333--342.}

\bibitem{Sul}
D.~Sullivan,
\newblock{Quasiconformal homeomorphisms and dynamics, I, Solution of
	the Fatou-Julia problem on wandering domains,}
\newblock{{\em Annals of Math.\/} {\bf 122} (1985), 401--418.}

\bibitem{Tat}
J.~Tate,
\newblock{Variation of the canonical height of a point
        depending on a parameter,}
\newblock{{\em Amer. J. Math.} {\bf 105} (1983), 287--294.}
\end{thebibliography}

\end{document}